\documentclass[oneside,a4paper,12pt,reqno]{amsart}
\usepackage{upref,amsmath,amsfonts,amsxtra,bbm, amssymb}
\newcommand{\Z}{\mathbb{Z}}
\newcommand{\T}{\mathbb{T}}
\newcommand{\R}{\mathbb{R}}
\newtheorem{theorem}{Theorem}[section]
\newtheorem{lemma}[theorem]{Lemma}
\newtheorem{prop}[theorem]{Proposition}
\newtheorem{cor}[theorem]{Corollary}
\begin{document}
\title{Periodicity and decidability of tilings of $\Z^2$}
\author{Siddhartha Bhattacharya}
\address{School of Mathematics, Tata Institute of Fundamental
  Research, Mumbai 400005, India}
\email{sid004@gmail.com}
\begin{abstract}
We prove that any finite set $F\subset \Z^2$ that tiles $\Z^2$ by
translations also admits a periodic tiling. As a consequence, 
the problem whether a given finite set $F$ tiles $\Z^2$ is decidable.
\end{abstract}
\maketitle

\section{Introduction}

Let $G$ be an abelian group and let $A, B$ and $C$ be subsets of $G$. We
will write $A\oplus B = C$ if every element of $C$ can be
uniquely expressed as $a + b$, with $a\in A$ and $b\in B$. We will
denote the set $\{ a + b : a\in A, b\in B \}$ by $A + B$.
If $d\ge 1$, and $F\subset \Z^d$ is a finite set
then a  {\it tiling \/} of $\Z^d$ by $F$ is a subset 
$C$ of $\Z^d$ satisfying $F\oplus C = \Z^d$. A tiling
$C$ is said to be {\it periodic\/} if there exists a finite index subgroup 
$\Lambda\subset \Z^d$ such that
$C + \Lambda = C$. 

For a finite set $F\subset \Z^d$ we define $\overline{F}\subset \R^d$
by ${\overline{F}} = (0,1]^d + F$.
It is easy to see that $\overline{F}$ tiles $\R^d$  if and only if 
$F$ tiles $\Z^d$. Questions about translational
tilings of $\R^d$ by bounded measurable subsets can often be reduced
 to the discrete case mentioned above (\cite{LW2},\cite{T}), and 
there is a large body of literature on 
tilings of $\Z^d$ by
translations of finite subsets (see \cite{S}, and the references
therein). In the $d=1$ case,
it is known that any tiling of $\Z$ by a finite
set $F$ is periodic with the period bounded by
$2^{\mbox{diam}(F)}$ (see\cite{N}). In particular, in this case the tiling problem 
is decidable. For $d \ge 1$, the following conjecture was proposed by
several authors (\cite{LW1}, \cite{SS}) :

\medskip
{\bf Periodic tiling conjecture : }{\it If a finite set $F\subset \Z^d$
tiles $\Z^d$  then it also admits a periodic tiling. \/}

\medskip
This conjecture is known to hold if $|F|$ is prime or $|F| =
4$ (see \cite{SZ}). If $d = 2 $, it is also known to be true if 
$\overline{F}$ is
simply-connected (\cite{BN}, \cite{WL}). More
generally, this holds for tilings of $\R^2$ by topological disks
\cite{K}.

\medskip
In this paper we give an affirmative answer when $d
= 2$. 

\medskip
\begin{theorem}\label{main}
Let $F,C$ be subsets of $\Z^2$ such that $F$ is finite and $F\oplus C =
\Z^2$. Then there exists $C^{'}\subset \Z^2$ and a finite index subgroup 
$\Lambda^{'}\subset \Z^2$ such that $F\oplus C^{'} = \Z^2$ and $C^{'}
+ \Lambda^{'} = C^{'}$.
\end{theorem}

\medskip
A well known argument, originally due to H. Wang, shows that for a
large class of tilings including the tilings we are considering here, this kind
of periodicity results imply the decidability of the underlying 
tiling problem
( \cite{IK}, \cite{R}). Hence as a consequence of Theorem
\ref{main} we obtain the following :

\medskip
\begin{cor}
The problem whether a given finite set $F$ tiles $\Z^2$ by
translations is decidable.
\end{cor}

\medskip
Although Theorem \ref{main} is a purely combinatorial result about
subsets of $\Z^2$, the proof we present here is ergodic-theoretic.
The key ingredient is a theorem about partitions of 
ergodic $\Z^2$-spaces.

\begin{theorem}\label{erg}
Let $(X,\mu)$ be a probability space equipped with an ergodic action of
$\Z^2$,
 and let $A\subset X$ be a measurable set such that $X$ admits
a partition of the form $(g_1A, \ldots , g_nA)$. 
Then there exists a partition $(A_1,\ldots ,A_m)$ of
$A$ such that each $A_i$ is invariant under the action of 
some non-zero element of $\Z^2$.
\end{theorem}  

This paper is organized as follows : In section 2 we show that Theorem
\ref{erg} implies Theorem \ref{main}. In section 3 we prove a weaker 
$L^2$-version of Theorem \ref{erg}, where the components are
functions rather than sets.
In section 4 we
show that these functions generate factors where
$\Z^2$-acts by commuting quasi-unipotent toral automorphisms. It is 
known that all ergodic invariant measures for such actions are algebraic
(it is a special case of Ratner's theorem on invariant measures of
unipotent actions, but in our context this can
be directly deduced from Weyl's equidistribution theorem). We use
this rigidity property and a combinatorial argument
  to complete the proof of Theorem \ref{erg} in section 5.
  \\

 \section{Dynamical formulation}

\medskip
If $\Z^2$ acts on a set $X$, we will call a point
 $x\in X$ {\it periodic\/} if the subgroup 
$H_x = \{ g : g\cdot x = x\}\subset \Z^2$ has finite index.
The point $x$ will be called $1$-{\it periodic\/} if $H_x$ 
is non-trivial. For any function $f : X\rightarrow \mathbb{C}$, the function
$x\mapsto f(g^{-1}\cdot x)$ from $X$ to $\mathbb{C}$ will be denoted by $g\cdot f$. For 
$E\subset X$, $1_E$ will denote the indicator function of $A$.

\begin{prop}\label{mean}
Suppose $(E_1, \ldots , E_m)$ is a partition of $\Z^2$ and $h_1, \ldots
,h_m$ are elements of $\Z^2$ such that  $E_i + h_i = E_i$ for all $i$. If
no two $h_i$'s are proportional to each other 
 then each $E_i$ is a periodic set.
\end{prop}

{\bf Proof.} Let $\Lambda_{ij}$ be the subgroup generated by $h_i$ and
$h_j$. Since $h_i$ and $h_j$ are not proportional to each other,
$\Lambda_{ij}$ is a finite index subgroup of $\Z^2$. If $\Lambda$ 
is the intersection of all
$\Lambda_{ij}$'s then $\Lambda$ also has finite index. 
We fix $1\le i\le m$ and pick $j\ne i$. Any $g\in
\Lambda$ can be written as $ah_i + bh_j$ with $a,b\in \Z$. As $E_i +
h_i = E_i$, we see that
$$ (g + E_i)\cap E_j = bh_j + (E_i \cap (E_j - bh_j)).$$
As $E_j - h_j = E_j$ and $E_i\cap E_j = \emptyset$, it follows that 
$ (g + E_i)\cap E_j = \emptyset$ for all $g\in \Lambda$.
Hence $E_i + \Lambda$ is disjoint from $E_j$ for each $j\ne i$.
Since $\cup E_i = \Z^2$, this shows that $E_i = E_i + \Lambda$. 
 $\Box$ 

\medskip
\begin{prop}
Suppose $F, D$ and $E$ are subsets of $\Z^2$ such that $F$ is finite
and $F\oplus D = E$. If $D$ is 1-periodic and $E$ is periodic then
there exists a periodic set $D^{'}$ such that $F\oplus D^{'} = E$.
\end{prop}

{\bf Proof.}
We first consider the case when $E = \Z^2$ and $D + (m,0) = D$ 
for some $m\ge 1$. 
We choose an integer $k > \mbox{diam}(F)$. 
For any $s\in \Z$ and $r\ge 1$,
let $W(s,r)\subset
\Z^2$ denote the horizontal strip defined by
$$ W(s,r) = \{ (i,j) : s\le j < s + r \}.$$
As $D$ is invariant under the element $(m,0)$ it follows that there
are at most $2^{mk}$ possibilities for the set $D\cap W(s,k)$. Hence
there exists $s\in\Z$ and $l\ge 2k$ such that the set $D\cap W(s,k)$ 
is a translate of the set $D\cap W(s + l,k)$. Let $D^{'}$ be the
unique set satisfying 
$$D^{'}\cap W(s,l) =  D\cap W(s,l)\ \mbox{and}\ 
D^{'} + (0,l) = D^{'}.$$
 As $D^{'}$ is invariant
under both $(m,0)$ and $(0,l)$, it is periodic. 
 For any $g\in \Z^2$,
let $B(g, k)$ denote the set $\{ h\in \Z^2 : ||g - h||\le k\}$. We
note that
for all $g\in \Z^2$, the set $B(g,k)\cap D^{'} $ is a translate of
$B(h,k)\cap D$ for some $h\in\Z^2$. Since $F\oplus D = \Z^2$, this
shows that $F\oplus D^{'} = \Z^2$ as well. 

Now suppose $E = \Z^2$ and $h$ is an arbitrary element of $\Z^2$ such
that $D + h =  D$.
Let $T$ be an element of $\mbox{GL}(2,\Z)$ such that 
$T(h)$ is of the form $(m,0)$ for some $m\ge 1$. Then $T(F) \oplus
T(D) = \Z^2$. As $T(D) + (m,0) = T(D)$, by the
above argument there exists a periodic set $D_1$ such that
$T(F) \oplus D_1 = \Z^2$. Clearly $D^{'} = T^{-1}(D_1)$ is periodic 
and $F\oplus D^{'} = \Z^2$.

Next we consider the case when $E$ is of the form $\Lambda + s$ for
some finite index subgroup $\Lambda\subset \Z^2$ and 
$s\in \Z^2$. Replacing $F$ by $F - s$ if necessary, we may assume that
$E = \Lambda$. We choose $h\in \Z^2$ such that 
 $D + h = D$. If 
$T : \Lambda \rightarrow \Z^2$ is an isomorphism then $T(F) + T(D) = \Z^2$
and the desired conclusion follows  from the previous case.  

To prove the general case we note that any periodic set $E$ is of the
form $\Lambda \oplus S$, where $\Lambda\subset \Z^2$ is a finite index
subgroup, and $S$ is a finite set. We pick $h\in\Z^2$ such that 
$D + h = D$. Without loss of generality we may assume that $h\in\Lambda$.
The set $D$
admits a partition $\{ D_s : s\in S\}$ such that $F\oplus D_s =
\Lambda + s$ and $D_s + h = D_s$. 
Applying the previous case, for each $s\in S$ we choose
$D^{'}_{s}\subset \Z^2$ and a finite index subgroup $\Lambda_s$ 
such that $F\oplus D_{s}^{'} = \Lambda + s$ and $D^{'}_{s} + \Lambda_s
= D^{'}_{s}$. If $D^{'} = \cup D^{'}_{s}$
then $F\oplus D^{'} = E $ and $ D^{'} +
\cap\Lambda_s =  D^{'}$. $\Box$

\medskip
Suppose $F\subset \Z^2$ is a finite set. 
Our next result shows that to prove Theorem \ref{main} it is enough to 
find a tiling $C$  that has a weaker periodicity property.

\medskip
\begin{prop}\label{partial}
Suppose $F\subset \Z^2$ is finite and
$F\oplus C = \Z^2$. If there exists a partition 
 $(C_1, \ldots , C_m)$ of $C$ such that each $C_i$ is 1-periodic, then
there exists a periodic set 
$C^{'}\subset \Z^2$ such that $F\oplus C^{'} = \Z^2$.
\end{prop}

{\bf Proof.} 
We start with a partition $( C_1, \ldots , C_m)$
for which $m$ is minimal. We choose
non-zero elements 
$h_1, \ldots, h_m$ in $\Z^2$ such that $C_i + h_i = C_i$ for all $i$.
Suppose  $h_i$ and $h_j$ are proportional for some $i\ne j$. Let 
$h$ be an element of $\Z^2$ such that $h = rh_i = sh_j$ for some 
non-zero $r,s\in \Z$. Then  $C_i + h = C_i$ and $C_j + h = C_j$. Hence 
if we remove $C_i$ and $C_j$ and add $C_i\cup
C_j$ then the resulting partition also satisfies the condition mentioned
above. This contradicts the minimality of $m$ and shows 
that no two $h_i$'s are proportional to each other. 

For $1\le i\le m$ we define $E_i = F\oplus C_i$. Since 
$(C_1, \ldots , C_m)$ is a partition of $C$ and $F\oplus C = \Z^2$, it
is easy to see that $(E_1, \ldots , E_m)$ is a partition of $\Z^2$. 
Applying Proposition \ref{mean} we deduce that each $E_i$ is a
periodic set. By the previous proposition there exists periodic sets 
$C_{1}^{'},\ldots , C_{m}^{'}$ such that $F\oplus C_{i}^{'} = E_i$ for
all $i$.
Since $(E_1, \ldots , E_m)$ is a partition of $\Z^2$, the sets 
 $F\oplus C_{i}^{'}  =  E_i$ and $F\oplus C_{j}^{'} = E_j$ are
disjoint whenever $i\ne j$. Hence the sets $C_{1}^{'}, \ldots ,
C_{m}^{'}$ are pairwise disjoint. We define $C^{'} = \cup C_{i}^{'}$ and
note that 
$$ F\oplus C^{'} = \cup (F\oplus C_{i}^{'}) = \cup (F\oplus C_{i}) =
\Z^2.$$
As $C^{'}$ is periodic, this proves the proposition.
$\Box$

\bigskip
Now we show that Theorem \ref{erg} implies Theorem \ref{main}. 
Let $F = \{g_1, \ldots ,g_n\}$ be a finite subset of $\Z^2$.  
We equip $\{ 0,1\}^{\Z^2}$ with the product topology 
and define $X(F)\subset \{ 0,1\}^{\Z^2}$ by
$$ X(F) = \{ 1_C : F \oplus C = \Z^2\}.$$

It is easy to see that $x\in X(F)$ if and only if for each $g\in
\Z^2$ there exists exactly one $g^{'}\in  g - F$ such that
$x(g^{'}) = 1$. This shows that $X(F)$ is a closed 
subset of the compact space $\{ 0,1\}^{\Z^2}$. Moreover, $X(F)$ is
invariant under the shift action of $\Z^2$ defined by
$(g\cdot x)(h) = x(h -g)$.  Since $\Z^2$ is amenable,
$X(F)$ admits ergodic shift invariant measures. 
We pick one such measure $\mu$ 
and define 
$$A = \{ x\in X(F) : x(0,0) = 1\}.$$
 Clearly 
$(g_1A, \ldots ,g_nA)$ is a partition of $X(F)$. By
Theorem \ref{erg}, the set $A$ admits a partition $(A_1,\ldots , A_m)$ such that
each $A_i$ is 1-periodic.
For any $x\in X(F)$ we define sets $C^{x},C_{1}^{x}, \ldots ,
C_{m}^{x}\subset \Z^2$ by 
$$ C^{x} = \{ g : g\cdot x \in A\},\ \ C_{i}^{x} = 
\{  g : g\cdot x \in A_i\}.$$

For each $i$, we choose  $h_i$ such that 
$h_i A_i = A_i$. Then $C_{i}^{x} + h_i = C_{i}^{x}$
 for $\mu$-almost all $x$.
As $\sum_{g\in F} 1_{gA} = 1$ and $\sum 1_{A_i} = 1_A$, we also see
that for $\mu$-almost all $x$, $F\oplus C^{x} = \Z^2$, and 
the sets $ C_{1}^{x}, \ldots ,
C_{m}^{x}$  form a partition of $C^{x}$.
In particular, there exists $x\in X$ such that  $F\oplus C^{x} = \Z^2$,
$C_{i}^{x} + h_i = C_{i}^{x}$ for all $i$, and $(C_{1}^{x}, \ldots ,
C_{m}^{x})$ is a partition of $C^{x}$.  
Now, Theorem \ref{main} follows from the previous proposition.

\medskip
\section{Spectral decomposition}

\medskip
Let $R$ denote the ring ${\mathbb
  Z}[u_{1}^{\pm}, u_{2}^{\pm}]$, the ring of
Laurent polynomials in two commuting variables with integer
coefficients. For $g = (n_1 ,n_2)\in {\mathbb Z}^2$, we
will denote the monomial
$u_{1}^{n_{1}}u_{2}^{n_{2}}$ by $u^{g}$. Every
element of $R$ can be expressed as $\sum c_{g} u^{g}$,
with $c_{g}\in \Z$ and $c_{g} = 0$ for all but
finitely many $g$. For $p = \sum c_{g} u^{g}$ and $r\ge 1$, we define 
$p^{(r)}\in R$ by
$ p^{(r)} = \sum c_{g} u^{rg}$. 
It is easy to see that for any $r$ the map $p\mapsto p^{(r)}$ is a
ring endomorphism. 
For 
$p = \sum c_gu^{g}\in
R$, we define $||p|| = \sum |c_g|$.  

Suppose $(X,\mu)$ is a probability space equipped with a
measure preserving action of $\Z^2$. 
For $p = \sum_{g} c_{g}u^{g}$ and $f\in
L^{2}(X)$ we define $p\cdot f\in L^{2}(X)$ by 
$ p\cdot f(x) = \sum_{{\bf g}} c_{g} f(g^{-1}\cdot x)$.  
It is easy to see that $L^{2}(X)$ becomes a $R$-module with respect
to this operation. For  $f\in L^{2}(X)$ we define an ideal 
$\mbox{Ann}(f)\subset R$ by
$$ \mbox{Ann}(f) = \{ p \in R : p\cdot f = 0\}.$$

Let $A\subset X$ be a measurable set such that $X$ admits a partition
of the form $(g_1A, \ldots , g_nA)$. Then $p_0\cdot 1_A = 1_X$ where
 $ p_0 = \sum_{i}u^{g_i}$. Hence $(u^{g}p_0 - p_0)\cdot 1_A = 0$  
for any $g\in\Z^2$. In particular, $\mbox{Ann}(1_A)$ is non-zero.
In this section we will show that any integer-valued function
$f\in L^{\infty}(X)$ with 
$\mbox{Ann}(f) \ne
0$ can be decomposed 
into a sum of bounded 1-periodic functions. We begin with a simple
observation.

\medskip
\begin{prop}\label{ann}
Let $f\in L^{\infty}(X)$ be an integer-valued function and
let $p$ be an element of $\mbox{Ann}(f)$. If 
$n_p$ denotes the product of all primes less
than $||p||||f||_{\infty} + 1$, then $p^{(r)}\in
\mbox{Ann}(f) $ whenever $r = 1 $ ( mod $n_p$).
\end{prop}

{\bf Proof.} Let $q$ be a prime number greater than $||p||||f||_{\infty}$.
It is easy to see that if $p = \sum c_g u^{g}$ then
$ p^{q} =  (\sum c_g u^{g})^{q} = p^{(q)} \ (\mbox{mod}\ q)$.
Hence $ p^{(q)} \cdot f = 0\  (\mbox{mod}\ q)$. Since 
$$|( p^{(q)}\cdot f)(x)| \le ||p||||f||_{\infty}$$ 
for almost all $x$, we deduce that $ p^{(q)}\cdot f = 0$, i.e., 
$ p^{(q)} \in \mbox{Ann}(1_A)$. 
Now let $r$ be any integer of the
form $j n_p + 1$. Then each prime factor of $r$ is
 greater than $||p||||f||_{\infty}$. By applying the
above argument repeatedly we conclude that $p^{(r)}\in \mbox{Ann}(f)$.
$\Box$ 

\medskip
For $g = (m,n)$ let $\chi_g$ denote the
character defined by  
$\chi_g (w_1,w_2) = w_{1}^{m}w_{2}^{n}.$
For any $s = \sum c_g u^{g}$ in $R$ we define $s_{*} :\T^2\rightarrow
\mathbb{C}$ by $s_{*}(w) = \sum c_g\chi_{g}(w)$. 
For any ideal $I\subset R$ we define 
$$V(I) = \bigcap_{s\in I} \{ w \in \T^2: s_{*} w = 0\}.$$
If $\Delta\subset\Z^2 - \{ 0\}$ is finite then $K(\Delta)$ will
denote the set $\cup_{g\in\Delta} \mbox{Ker}(\chi_g)$.

\medskip
\begin{lemma}\label{variety}
Let $(X,\mu)$ be a probability space together with a measure
preserving action of $\Z^2$, and let $f$ be an integer valued function
in $L^{\infty}(X)$ with $\mbox{Ann}(f) \ne 0$. Then there exists a
finite set $\Delta\subset \Z^2 -\{ 0\}$ such that
no two elements of $\Delta$ are proportional to each other
and $V(\mbox{Ann}(f)) \subset K(\Delta)$.
\end{lemma} 

{\bf Proof.} We choose a non-zero $p = \sum c_g u^{g}\in\mbox{Ann}(f)$. Since
$u^{g}p \in\mbox{Ann}(f)$ for all $g\in\Z^2$, without loss of
generality we may assume that $c_0 \ne 0$. 
We define $\mbox{Supp}(p) = \{ g \in \Z^2 : c_g\ne 0\}$ and 
$E = n_p (\mbox{Supp}(p) - \{ 0\})$. Let $w$ be an
element of $V(\mbox{Ann}(f))$. From the previous proposition it follows that
 $p^{(r)}_{*}(w) = 0$ whenever $r$ is of the form  $jn_p + 1 $.
Hence for all $j\ge 1$,
$$ \sum_{g\ne 0}c_g\chi_{g}^{jn_p + 1}(w) = p^{(jn_p + 1)}_{*}(w) - c_0 =
-c_0.$$
If we put $j = 1,2,\ldots , N$ and take the average, this shows
that
$$ \sum_{g\ne 0}c_g\chi_{g}(w)({\frac{1}{N}}
\sum_{j=1}^{N}\chi_g(w)^{jn_p}) = -c_0.$$
Hence there exists $g\in \mbox{Supp}(p) - \{0\}$ such that 
${\frac{1}{N}}
\sum_{j=1}^{N}\chi_g(w^{n_p})^{j}$ does not go to zero as $N\mapsto
\infty$. Clearly for any such $g$, $ \chi_g(w^{n_p}) = 1$. Since 
$ \chi_g(w^{n_p}) = \chi_{n_p g}(w)$, this proves that $w\in K(E)$. 
Hence $V(\mbox{Ann}(f))$ is contained in $K(E)$, and 
we have shown that there exist finite sets that satisfy the
second condition. Let $\Delta\subset\Z^2$ be such a set with 
least number of elements. Suppose there exist $g\ne h\in \Delta$ such that 
$qg = rh = g_1$ for some non-zero $q,r\in\Z$. Then $\mbox{Ker}(\chi_{g_1})$
contains both $\mbox{Ker}(\chi_{g})$  and $\mbox{Ker}(\chi_{h})$. Let
$\Delta^{'}$ denote the set obtained from $\Delta$ by removing $g$ and
$h$, and adding $g_1$. As $K(\Delta^{'})\supset K(\Delta)$, the set 
$\Delta^{'}$ also satisfies the second condition. This contradicts the
minimality of $|\Delta |$ and shows that no two elements of $\Delta$
are proportional to each other. $\Box$

\medskip
For $g\in \Z^2$, we define $H_g\subset L^{2}(X)$ by 
 $H_g = \{ h : u^{g}\cdot f = h\}$. For $f\in L^{2}(X)$, $f^{g}$
will denote the projection of $f$ onto $H_g$. By the pointwise ergodic
theorem
$$ {\frac{1}{n}} \sum_{i=0}^{n - 1} (u^{ig}\cdot f)(x) \mapsto f^{g}(x) \
  \mbox{as}\ n\mapsto \infty$$
for almost all $x$.

\medskip
\begin{theorem}\label{l2}
Let $(X,\mu)$ be a probability space together with a measure
preserving action of $\Z^2$, and let $f$ be an integer valued function
in $L^{\infty}(X)$ with $\mbox{Ann}(f) \ne \{ 0\}$. Then there exists a
finite set $\Delta\subset \Z^2 -\{ 0\}$ such that
\begin{enumerate}
\item{No two elements of $\Delta$ are proportional to each other}
\item{The function $f_0 = f - \sum_{g\in  \Delta}f^{g} $ is 
 periodic}
\item{For any $l\ge 1$ and $g\in \Delta$, the function $f^{lg} -
    f^{g}$ is periodic.}
\end{enumerate}
\end{theorem}

\medskip
{\bf Proof.} Let $H\subset L^2(X)$ be the closure of the linear
subspace generated by
the set $\{ u^{g}\cdot f : g\in\Z^2\}$. For any
probability measure $\nu$ on ${\mathbb T}^2$ we define a unitary
representation $\sigma_{\nu}$ of $\Z^2$ on $L^{2}(\T^2, \nu)$ by
$\sigma_{\nu}(g)(\phi)(x) = \chi_g(x)\phi(x)$.
We note that $f$ is a cyclic vector for the representation of 
$\Z^2$ on $H$. By the spectral
theorem for commuting unitary operators, there exists a probability
measure $\nu$ on ${\mathbb T}^{2}$ and a unitary isomorphism 
$\theta : L^{2} ({\mathbb T}^{2}, \nu) \rightarrow H$ such that 
$\theta \circ \sigma_{\nu}(g)(\phi) = u^{g}\cdot \theta(\phi)$ 
for all $g\in \Z^2$,
and $\theta(1_{{\mathbb T}^{2}}) = f$. Clearly for any $s\in
\mbox{Ann}(f)$ we have
$$ \theta(s_{*}) = \theta(s_{*}\cdot 1_{\T^2}) = s\cdot f = 0.$$
Hence $\nu$ is supported on the set $\{ w : s_{*}(w) = 0\}$. Since
this holds for all $s\in \mbox{Ann}(f)$, this shows that 
$\nu$ is supported on $V(\mbox{Ann}(f))$. Let $\Delta \subset \Z^2 -
\{ 0\}$ be as in the previous lemma. Then $\nu(K(\Delta))= 1$ and
$f = \theta(1_{K(\Delta)})$. For $g\in \Delta$ we define 
$H^{'}_{g} \subset L^{2}(\T^2,\nu)$ by 
$H^{'}_{g} = \{ \phi : \phi\circ\sigma_{\nu}(g) = f\}$. 
It is easy to see that $\phi\in H^{'}_{g}$
if and only if it vanishes outside the set $K_g = \mbox{Ker}(\chi_g)$.
As $\theta$ is a conjugacy it follows that $f^{g} =
\theta(1_{K_g})$. This implies that 
$$f  - \sum_{g\in
  \Delta}f^{g}  = \theta(1_{\cup K_g} - \sum_{g\in
  \Delta} 1_{K_g}).$$

Let $B$ denote the set of all $w\in \T^2$ that belongs to at least
two $K_g$'s. Since no two $g$'s are proportional to each other,
 $B$ is a finite set and each element of $B$ is a root of
unity. Hence there exists a finite index subgroup $\Lambda\subset
\Z^2$ such that $\chi_{g}(b) = b$ for all $b\in B$ and $g\in\Lambda$. Since the
function $1_{\cup K_g} - \sum_{g\in
  \Delta} 1_{K_g}$ is supported on $B$, it is invariant under 
$\sigma_{\nu}(\Lambda)$. 
As $\theta$ is a conjugacy, we conclude that 
the function $f_0 = f - \sum_{g\in
  \Delta}f^{g} $ is also $\Lambda$-invariant. 

To prove the last assertion we note that $f^{lg} - f^{g} = 
\theta(1_{K_lg} - 1_{K_g})$. We define $Q\subset \T^2$ by
$Q = \cup\{ K_{lg}\cap K_{h} : h\in \Delta - \{g\} \}.$
Since distinct elements of $\Delta$ are not proportional to each
other, $Q$ is a finite set consisting of roots of unity. Hence any
function supported on $Q$ is periodic. As $\nu(K(\Delta)) = 1$, it is
easy to see that there is a function $f^{*}$ such that the support of
$f^{*}$ is contained in $Q$ and $1_{K_lg} - 1_{K_g} = f^{*}$
a.e. $\nu$. Since $\theta $ is a conjugacy, we conclude that 
$f^{lg} - f^{g}$ is periodic. $\Box$\\

\section{Ergodic decomposition}

Throughout this section $(X,\mu)$ will denote a probability space
equipped with an ergodic action of $\Z^2$, and $A\subset X$ will
denote a measurable set such that $\mbox{Ann}(1_A)$ is non-trivial.
Suppose $\Lambda\subset \Z^2$ is a finite index subgroup, and
 ${\mathcal A}$ is the $\sigma$-algebra of all $\Lambda$-invariant sets. Since
the induced $\Z^2$-action on $(X, {\mathcal A}, \mu)$ is ergodic and
it factors through an action of the finite group $\Z^2/\Lambda$, it
follows that ${\mathcal A}$ is finite. Let  $\{ X_s : s\in
S(\Lambda)\}$ be the collection of atoms of ${\mathcal A}$. For each 
$s\in S(\Lambda)$ we define a measure $\mu_s$ by $\mu_s(E) =
\mu(X_s\cap E) / \mu(X_s)$. Elements of  $\{ X_s : s\in
S(\Lambda)\}$ will be called $\Lambda$-ergodic components of $X$, and 
elements of $\{ \mu_s : s\in S(\Lambda)\}$ will be called
$\Lambda$-ergodic components of $\mu$.

\medskip  
Let $T$ be a finite dimensional torus and let $\alpha$ be an action of
$\Z^2$
on $T$ by continuous automorphisms. The action $\alpha$ is called 
 {\it unipotent \/} if for each $g$
there exists $k\ge 1$ such that $(\alpha(g) - I)^k = 0$. We will use
the following result about invariant measures of such actions.
\begin{prop}\label{weyl}
Let $\alpha$ be as above and let $\nu$ be an ergodic
$\alpha$-invariant measure. Then there exists a finite index subgroup
$\Lambda$ such that each $\Lambda$-ergodic component of $\nu$ 
is the Haar measure on a translate of a subtorus.
\end{prop}

 Since every $\alpha$-invariant
ergodic measure arises from an orbit, and for each $x\in T$ the map 
$g\mapsto \alpha(g)(x)$ is a polynomial map, one can  deduce this
result from the multi-dimensional version of Weyl's equidistribution
theorem (see \cite{T2} for details). 

For $\phi : X\rightarrow \R$ we define $\overline{\phi} : X\rightarrow \T$ by 
$\overline{\phi} (x) = e^{2\pi i \phi(x)}$. 
Let $f_0$ and 
$\{ f^g : g\in\Delta \}$ denote the functions
associated with the decomposition of $1_A$ described in Theorem \ref{l2}. 
 For any $g\in\Delta$,  $\mbox{Ann}({\overline{f^g}})$ is contained in 
$\mbox{Ann}(f^g)$. In particular, $u^{g} - 1 \in
\mbox{Ann}({\overline{f^g}})$. 
The next lemma shows
that $\mbox{Ann}(\overline{f^g})$ also contains elements that are 
not divisible by 
 $u^{g} - 1 $.

\begin{lemma}\label{torus}
Let $g\in\Delta$ and $f^g$ be as above. Then
there exists $h\in\Z^2$ such that $h$ is not proportional to $g$ and 
 $(u^{h} - 1)^k \in \mbox{Ann}(\overline{f^g})$ for some $k\ge 1$.
\end{lemma}
 
{\bf Proof.} We choose any $g^{'}$ that is not proportional to
$g$. Let $\Lambda_1$ be the finite index subgroup generated by
$g$ and $g^{'}$, and let $\Lambda_2$  be a finite index subgroup
that leaves $f_0$ invariant. We pick $l\ge 1$ such that
$l\Z^2$ is contained in $\Lambda_1\cap \Lambda_2$, and define $q\in
R$ by
$$q = \prod_{g_1\in\Delta -\{g\}}(u^{lg_1} -1).$$
We note that  for any $g_1\in\Delta -\{g\}$, the polynomial $u^{lg_1} -1$
  annihilates both $f^{g_1}$ and $f_0$. This implies that 
$$ q\cdot f^g = q\cdot (f_0 + \sum_{g\in\Delta}f^g) = q\cdot
1_A.$$
Hence $q\cdot f^g(x)\in \Z$ for almost all $x$, i.e., 
 $q \in \mbox{Ann}(\overline{f^g})$.
We fix any $g_1\in \Delta -\{g\}$. 
By our choice of $l$, $lg_1 = rg + r^{'}g^{'}$ for some 
$r,r^{'}\in \Z$. This shows that
 $$u^{lg_1} - 1 = u^{r^{'}g^{'}} (u^{rg} -1) + 
( u^{r^{'}g^{'}} - 1).$$
 Since  $(u^{rg} -1)\cdot {\overline{f^g}} = 0$, it follows
  that $( u^{r^{'}g^{'}} - 1)\cdot {\overline{f^g}} = 
( u^{lg_1} - 1)\cdot {\overline{f^g}}$. 
We choose a number $N $
    such that $l$ and all the $r^{'}$'s corresponding to different
    values of $g_1$ divides $N$. As $q\in \mbox{Ann}({\overline{f^g}})$,
 it follows that $(u^{Ng^{'}} - 1)^{k} \cdot
     {\overline {f^g}} = 0$ where $k = |\Delta | -1 $. Hence  
$h = Ng^{'}$ satisfies the given condition. $\Box$  

\medskip
\begin{lemma}
Let $(X,\mu)$ and $A\subset X$ be as before.
Then there exists a finite index subgroup $\Lambda\subset\Z^2$ such
that for each $\Lambda$-ergodic component $\mu_s$ of $\mu$ 
the function $f_0$ is a constant a.e. $\mu_s$,
and for all $g\in \Delta$ either ${\overline{f^g}}$ is a constant a.e. $\mu_s$ 
or ${\overline{f^g}}_{*}(\mu_s)$ is the Haar measure on $\T$.
\end{lemma}

{\bf Proof.}
Let $Y$ denote the compact abelian group $\T^{\Z^2}$ and let 
$\alpha$ denote the shift action of $\Z^2$ on $Y$. For
$g\in\Delta$, let $h$ and $k$ be as in the previous lemma. We define 
a closed subgroup $H_g\subset Y$ by
$$ H_g = \{ y : (\alpha(g) - I)y = (\alpha(h) - I)^{k}y = 0\}.$$
Since $g$ and $h$ are not proportional to each other, $H_g$ is a finite
dimensional torus. If $\Lambda_g$ denote the subgroup generated by $g$ 
and $h$ then it is easy to see that $\Lambda_g$ acts on $H_g$ by
unipotent automorphisms. For $g\in \Delta$ we define 
$\pi^{g} : X\rightarrow Y$ by 
$$ \pi^{g}(x)(g_1) = {\overline{f^g}}(g_{1}^{-1}\cdot x).$$
From the previous lemma we deduce that $\pi^{g}_{*}(\mu)$ is supported 
on $H_g$. We define $ H = \prod H_g$ and $\pi = \prod
\pi^{g}$. Clearly $\pi : X\rightarrow H$ is an equivariant map. Hence
$\pi_{*}(\mu)$ is ergodic with respect to $\alpha$. Let $\Lambda_0$ be
a finite index subgroup that leaves $f_0$ invariant, and let $\Lambda_1$ 
be the intersection of $\Lambda_0$ and all $\Lambda_g$'s. Then
$\Lambda_1$ is a finite index subgroup and each element of $\Lambda_1$
acts on $H$ by unipotent automorphisms. By Proposition \ref{weyl} 
there exists a finite index subgroup $\Lambda\subset \Lambda_1$ such that
 any $\Lambda$-ergodic component $\pi_{*}(\mu)_{s_1}$ of $\pi_{*}(\mu)$
is supported on a set $L_{s_1}$ that is a translate of a subtorus, and
$\pi_{*}(\mu)_{s_1}$ is the Haar measure on $L_{s_1}$.  
Let $\mu_{s}$ be a $\Lambda$-ergodic component of $\mu$. 
As $\{ \pi^{-1}(L_{s_1}) :s_1\in S(\Lambda)\}$ is a partition of $X$ by
$\Lambda$-invariant sets, we deduce that $X_{s}\subset
\pi^{-1}(L_{s_1})$ for some $s_1$. 
Hence $\pi_{*}(\mu_{s})$ is an
ergodic $\Lambda$-invariant measure on $L_{s_1}$. Since $\mu_{s}$ is
absolutely continuous with respect to the restriction of $\mu$ to 
$\pi^{-1}(L_{s_1})$, it follows that
$\pi_{*}(\mu_{s})$ is  absolutely continuous with respect to 
$\pi_{*}(\mu)_{s_1}$, the Haar measure on $L_{s_1}$. On the other hand, any
proper subset of $L_{s_1}$ that is a translate of a subtorus has zero Haar
measure. By Proposition 
\ref{weyl}, $\pi_{*}(\mu_{s})$ is the Haar measure on $L_{s_1}$.
Clearly, $f^{g}_{*}(\mu_{s})$ is the projection of
$\pi_{*}(\mu_{s})$ to the $(0,0)$'th co-ordinate of $H_g$.
As $L_{s_1}$ is connected, we conclude that the measure
$f^{g}_{*}(\mu_{s})$ is either supported on a point or is the Haar
measure on $\T$. $\Box$
  
\medskip 
We will use the previous lemma to prove a weaker version of Theorem
\ref{erg} that holds for any $A$ with ${\mbox{Ann}}(1_A)\ne \{ 0\}$. 
\begin{theorem}\label{half}
Let $(X,\mu)$ be a probability space equipped with an ergodic
$\Z^2$-action, and let $A\subset X$ be a measurable subset such that
${\mbox{Ann}}(1_A)\ne \{ 0\}$. Then there exists a finite index
subgroup $\Lambda\subset \Z^2$ such that for each $\Lambda$-ergodic
component $X_s\subset X$, either $X_s\cap A$ is 1-periodic or
 $\mu(X_s\cap A) = \mu(X_s)/2$.
\end{theorem}
{\bf Proof.} 
 Let $\Lambda$ be as described above. We fix $s\in
S(\Lambda)$. Suppose $\overline{f^g}$ is constant with respect to
$\mu_s$ for all $g\in \Delta$. We claim that in this case $X_s\cap A$
is 1-periodic. We first consider the case when all these
constants are different from 1. Since the function $t\mapsto e^{2\pi i
  t}$ is injective on $(0,1)$, it follows that each $f^g$ is also a
constant. Since $f_0$ is also a constant on $X_s$, the same holds for 
$1_A = f_0 + \sum f^{g}$. Hence either $X_s\cap A = \emptyset$ or 
$X_s\cap A = X_s$. Since $X_s$ is a periodic set, 
so is $X_s\cap A$.

Suppose  $\overline{f^g}$ is constant with respect to
$\mu_s$ for all $g\in \Delta$ and $\overline{f^g} = 1$ on $X_s$ for
some $g$.
Then $f^{g}(x)\in \{ 0,1\}$ for $\mu_s$-almost all $x$. We choose
$l\ge 1$ such that $lg\in\Lambda$. By the pointwise ergodic theorem,
for all $h\in\Z^2$ we have 
$$ f^{h}(x) = \lim_{n\mapsto\infty}{\frac{1}{n}}\sum_{i=0}^{n-1}
( u^{ih}\cdot 1_A)(x)$$
$\mu_s$-almost everywhere. This shows that for $\mu_s$-almost all $x$, 
$$ f^{g}(x) = {\frac{1}{l}}\sum_{i=0}^{l-1}f^{lg}(-ig\cdot x).$$
As $f^{lg}(y) \in [0,1]$ for all $y$, the right hand side is equal to 1 if and
only if each term is equal to 1. This implies that $f^{lg}(x) = 1$ whenever
$f^{g}(x) = 1$. Similarly,  $f^{lg}(x) = 0$ 
whenever $f^{g}(x) = 1$. Hence $f^{lg}(X_s)\subset \{0,1\}$. 
Since $\mu_s$ is $\Lambda$-invariant we deduce that 
$f^{lg}$ restricted to $X_s$ is the projection of $1_{X_s\cap A}\in
L^{2}(X_s, \mu_s)$ onto the  $\sigma$-algebra
 $lg$-invariant functions in $L^{2}(X_s, \mu_s)$.
Hence $\int_E f^{lg} d\mu_s = \mu_s(A\cap E)$ for any $lg$-invariant
subset $E$ of $X_s$. Putting $E = \{ x \in X_s : f^{g}(x) = 1\}$ we
see that 
$ \mu_{s}(E) = \int_{E}f^{lg} d\mu_s  = \mu_s(E\cap A)$.
Hence  $A\supset E$. We also have
$$\mu_{s}(E) = \int_{X_s} f^{lg} d\mu_s = \mu_s(A).$$
Hence  $E = A\cap X_s$. Since $lg\cdot E = E$, 
it follows that $lg\cdot (A\cap X_s) = A\cap X_s$. 
This proves the claim.

Now suppose ${\overline{f^g}}_{*}(\mu_s)$ is the Haar measure on $\T$
for some $g\in \Delta$. Then $f^g_{*}(\mu_s)$ is the Haar measure on
$[0,1]$. As in the previous case, we choose $l\ge 1$ such that 
$lg\in \Lambda$.
By Theorem \ref{l2}, $f^{lg} - f^{g}$ is
invariant under a finite index subgroup $\Lambda_1$. We define
$\Lambda_2 = \Lambda\cap \Lambda_1$. Let $\{ Y_1, \ldots, Y_r\}$
denote the $\Lambda_2$-ergodic components of $(X_s, \mu_s)$.  For $i = 1,
\ldots , r$; there exists a constant $c_i$ such that $f^{lg} - f^{g} =
c_i$ a.e. on $Y_i$. For each $i$, we define a measure $\nu_i$ on $X_s$ by 
$\nu_i(E) = \mu_s(E\cap Y_i)/\mu_s(Y_i)$. Clearly $\pi_{*}(\nu_i)$ is
$\Lambda_2$-invariant and ergodic with respect to $\Lambda_2$. As 
$\nu_i$ is absolutely continuous with respect to $\mu_s$, we deduce
that $\pi_{*}(\nu_i)$ is absolutely continuous with respect to 
$\pi_{*}(\mu_s)$. Since $\pi_{*}(\mu_s)$ is the Haar measure on a
translate of a subtorus, and the $\Lambda_2$-action on $H$ is
unipotent, by Proposition \ref{weyl}, $\pi_{*}(\nu_i) =
\pi_{*}(\mu_s)$. Hence $f^{g}_{*}(\nu_i)$ is also the Lebesgue measure
on $[0,1]$ and $f^{lg}_{*}(\nu_i)$ is the Lebesgue measure
on $[c_i, c_i + 1]$. On the other hand,
$f^{lg}_{*}(\nu_i)$  is a probability neasure on $[0,1]$. Combining
these these two observations we deduce that $c_i = 0$ for all $i$,
i.e., $f^{lg} = f^{g}$ on $X_s$.
We have already observed that $ \int f^{lg}d\mu_s = \mu_s(A)$. Hence
$\mu_s(A) =  \int f^{g}d\mu_s = {\frac{1}{2}}$. 
$\Box$ 

\medskip  
\section{Proof of the main theorem}

In this section  $A\subset X$ will
denote a measurable set such that $X$ admits a partition of the form
$(g_1A, \ldots , g_nA)$ with $\{ g_1, \ldots, g_n\} = F \subset \Z^2$.
This implies that ${\mbox{Ann}}(1_A)$ is
non-trivial. As before, $f_0$ and 
$\{ f^g : g\in\Delta \}$ will denote the functions
associated with the decomposition of $1_A$ described in Theorem \ref{l2}. 

Recall that a set $E\subset X$ is periodic if $\{g : gE = E\}$ 
is a finite index subgroup of $\Z^2$, and it is 1-periodic if
$gE = E$ for some non-zero $g\in \Z^2$. We will call a set $E$
{\it weakly periodic \/} if it can be expressed as finite disjoint union of 
1-periodic sets. 

\begin{lemma}
If $E\subset X$ is periodic and $B\subset E$ is weakly periodic then
$E - B$ is weakly periodic. 
\end{lemma}

{\bf Proof.} Let $m$ be the smallest positive integer such that there
exist disjoint 1-periodic sets $B_1, \ldots , B_m$ with 
$\cup B_i = B$. For each $i$ we choose a non-zero $g_i\in\Z^2$ such that $g_iB_i
= B_i$. Suppose there exist $i\ne j$ such that $g_i$ and $g_j$ are
proportional to  each other. Then the set $B_i\cup B_j$ is also
1-periodic. Hence if we remove $B_i$ and $B_j$ and add $B_i\cup B_j$
then the smaller collection also satisfies the above condition. Since
this contradicts the minimality of $m$ we deduce that no two $g_i$'s
are proportional to each other.

Let $\Lambda_{ij}$ be the subgroup generated by $g_i$ and
$g_j$, and let $\Lambda_0 = \{ g : gE = E \}$. 
Let $\Lambda$ denote the intersection
of $\Lambda_0$ and all $\Lambda_{ij}$'s. Clearly $\Lambda$ has finite
index. Let $(Q_1, \ldots , Q_r)$ be the partition
of $X$ by $\Lambda$-ergodic components. 
Since $E$ is invariant under 
$\Lambda_0\supset \Lambda$ and the restriction of the 
$\Lambda$-action to each $Q_i$ is ergodic, it follows that 
$E = \cup\{ Q_i : Q_i \subset E\}$. This implies that 
$$ E - B =   \cup \{ Q_i\cap B^c : Q_i \subset E\}.$$
So to prove $E - B$ is weakly periodic it is enough to show that 
for each $i$ the set $Q_i\cap B^c$ is 1-periodic. As $Q_i$ is 
$\Lambda$-invariant, this is trivially true if $B$ does not intersect
$Q_i$. Suppose $B$ intersects $Q_i$. Then $Q_i \cap B_j $ has positive
measure for some j. For $l\ne
j$, $B_l$ is invariant under $g_l$ and $B_j\cap B_l =
\emptyset$. Hence $ng_lB_j\cap B_l = \emptyset$ for all $n\in
\Z$. As $B_j$ is invariant under $g_j$, the set 
$Y = \cup \{ g B_j : g\in \Lambda_{lj}\}$ is disjoint from $B_l$. Since
$Y$ is $\Lambda_{ij}$-invariant and intersects $Q_i$ non-trivially, 
it contains $Q_i$. Hence $B_l\cap Q_i = \emptyset$. 
Since $l\ne j$ is arbitrary, this shows that $Q_i\cap B = Q_i\cap
B_j$, i.e., $Q_i\cap B^c = Q_i\cap
B_{j}^{c}$. As $Q_i$ is periodic and $B_j$  is 1-periodic, this
completes the proof. $\Box$

\medskip
We note that if Theorem \ref{erg} is true for $A$ then it
is also true for any translate of $A$. Replacing $F$ by a translate 
if necessary, we may assume that  $0\in F$ and 
$$F  \subset \{ (m,n) : n > 0\} \cup \{ (m,n) : n = 0, m\ge 0\}.$$  
Hence for any finite set of non-zero elements of $F$
their sum is non-zero.

\medskip
{\bf Proof of Theorem \ref{erg}.} 
We define $K = - F$. Let $\Lambda\subset \Z^2$ be a
finite index subgroup satisfying the condition stated in Theorem
\ref{half}. Clearly $\Z^2$ acts transitively on $S(\Lambda)$, the set
of $\Lambda$-ergodic components of $X$.
 For each $s\in S(\Lambda)$, we define $A_s = X_s\cap A$. 
It is enough to show that $A_s$ is
weakly periodic for all $s$. As 
$\{ gA : g\in F \}$ is a partition of $X$, it follows that 
for each $s\in S(\Lambda)$, $\{ g^{-1}A_{g\cdot s} : g\in K\}$ is a partition
of $X_s$.  Suppose $A_s$ is not
weakly periodic for some $s$. Then $\mu(A_s) = \mu(X_s)/2$. 
By the previous lemma the set $X_s\cap A_{s}^{c}$ is also not weakly
periodic. On the other hand, as $0\in F$, it follows that 
 $$X_s\cap A_{s}^{c} = \cup\{ g^{-1}A_{g\cdot s} : g\in K - \{0\}\}.$$
 Hence $A_{g\cdot s}$ is not weakly periodic for some 
$g\in K - \{0\}$. As   $\mu(A_s) = \mu(A_{g\cdot s}) = \mu(X_s)/2$,
this shows that $g^{-1}A_{g\cdot s} = X_s\cap A_s^{c}$, i.e.,
$A_{g\cdot s} = g( X_s -  A_s)$. We replace $A_s$ by $A_{g\cdot s}$
and apply this argument again to obtain $h\in K -\{ 0\}$ such that 
$A_{(g + h)\cdot s}$ is not weakly periodic and  
$A_{(g + h)\cdot s} = h( X_{g\cdot s} -  A_{g\cdot s})$. As 
$X_{g\cdot s} = gX_s$ and $A_{g\cdot s} = g( X_s -  A_s)$, it follows
that $A_{(g + h)\cdot s} = (g + h)\cdot A_s$. As both $g$ and $h$ are
non-zero,  $g+h$ is also non-zero.
We define $g_1 = g, h_1 = h$ and $s_1 = (g + h)\cdot s$. 
By repeated application of the above
argument for each $i\ge 1$ we 
obtain $ s_i, g_i$ and $h_i$ such that
 $A_{s_{i}} = (g_i + h_i)\cdot A_{s_{i-1}}$
and no $A_{s_i}$ is weakly periodic. 
Since $S(\Lambda)$ is finite there exists $i < j$ such that $A_{s_i} =
A_{s_j}$. Then $g_0A_{s_i} = A_{s_i}$ for some $g_0\in\Z^2$ that 
can be expressed as sum of non-zero elements of $K$. By our choice of
$F$, $g_0\ne 0$. Since this contradicts the fact that 
 no $A_{s_i}$ is weakly periodic, we conclude that 
$A_s$ is weakly periodic for all $s$. $\Box$

\end{document}